\documentclass{LMCS}

\def\doi{8 (3:06) 2012}
\lmcsheading%
{\doi}
{1--12}
{}
{}
{Mar.~22, 2011}
{Aug.~10, 2012}
{}

\usepackage{enumerate}
\usepackage{hyperref}

\usepackage[all,ps]{xy}

\usepackage{tikz}
\usetikzlibrary{arrows}

\input{finhilb.sty}

\begin{document}

\title[Finite dimensional Hilbert spaces are complete]{Finite dimensional Hilbert spaces are
  complete for dagger compact closed categories\rsuper*}

\author[P.~Selinger]{Peter Selinger}
\address{Department of Mathematics and Statistics\\
  Dalhousie University,
  Halifax, Nova Scotia, Canada}
\email{selinger@mathstat.dal.ca}
\thanks{Research supported by NSERC.}

\titlecomment{{\lsuper*}An extended abstract of this work has appeared in
  {\em Proceedings of the 5th International Workshop on Quantum Physics
  and Logic (QPL 2008)}, Electronic Notes in Theoretical
  Computer Science 270(1):113-119, 2011}

\begin{abstract}
  \noindent We show that an equation follows from the axioms of dagger
  compact closed categories if and only if it holds in finite
  dimensional Hilbert spaces.
\end{abstract}

\maketitle

\section{Introduction}

Hasegawa, Hofmann, and Plotkin recently showed that the category of
finite dimensional vector spaces over any fixed field $k$ of
characteristic $0$ is {\em complete} for traced symmetric monoidal
categories {\cite{HHP08}}. What this means is that an equation holds
in all traced symmetric monoidal categories if and only if it holds in
finite dimensional vector spaces. The authors also noted that it is a
direct corollary, via Joyal, Street, and Verity's ``Int''-construction
{\cite{JSV96}}, that finite dimensional vector spaces are also
complete for compact closed categories. The present paper makes two
contributions: (1) we simplify the proof of Hasegawa, Hofmann, and
Plotkin's result, and (2) we extend it to show that finite dimensional
{\em Hilbert} spaces are complete for {\em dagger} traced symmetric
monoidal categories (and hence for dagger compact closed categories).

The paper is organized as follows. We assume the reader knows the
definition of a dagger compact closed category {\cite{AC04,Sel05}},
and has at least an informal understanding of their graphical language
{\cite{Sel11}}. In Section~\ref{sec-statement}, we state the main
result without further ado.  Sections~\ref{sec-reductions} and
{\ref{sec-informal}} are devoted to an informal, but hopefully
comprehensible, explanation of the proof. For the reader who is
interested in details, full technical definitions and proofs
(including a formal definition of the graphical language and
isomorphism of diagrams) appear in
Section~\ref{sec-technical}. Section~\ref{sec-generalizations}
discusses how to generalize the result to fields other than the
complex numbers, gives counterexamples for some possible
strengthenings of completeness, and concludes with an open problem.

\section{Statement of the main result}\label{sec-statement}

For a definition of dagger compact closed categories, their term
language, and their graphical language, see
{\cite{AC04,Sel05,Sel11}}. We also use the concept of a {\em dagger
  traced monoidal category} {\cite{Sel11}}, which is a dagger
symmetric monoidal category with a trace operation {\cite{JSV96}}
satisfying $\Tr^X_{U,V}(f)\da = \Tr^X_{V,U}(f\da)$. We note that every
dagger compact closed category is also dagger traced monoidal;
conversely, by Joyal, Street, and Verity's ``Int'' construction, every
dagger traced monoidal category can be fully embedded in a dagger
compact closed category.

We will make use of the soundness and completeness of the graphical
representation, specifically of the following result:

\begin{thm}[\cite{Sel05}]\label{thm-graphical}
  A well-typed equation between morphisms in the language of dagger
  compact closed categories follows from the axioms of dagger
  compact closed categories if and only if it holds, up to graph
  isomorphism, in the graphical language.
\end{thm}

An analogous result also holds for dagger traced monoidal categories
{\cite[Thm.~7.12]{Sel11}}.
The goal of this paper is to prove the following:

\begin{thm}\label{thm-main}
  Let $M,N:A\ii B$ be two terms in the language of dagger compact
  closed categories. Suppose that $\semm{M}=\semm{N}$ for every
  possible interpretation (of object variables as spaces and morphism
  variables as linear maps) in finite dimensional Hilbert spaces. Then
  $M=N$ holds in the graphical language (and therefore, holds in all
  dagger compact closed categories).
\end{thm}

\section{Reductions}\label{sec-reductions}

Before attempting to prove Theorem~\ref{thm-main}, we reduce the
statement to something simpler.  By arguments analogous to those of
Hasegawa, Hofmann, and Plotkin {\cite{HHP08}}, it suffices without loss
of generality to consider terms $M,N$ that satisfy some additional
conditions. The additional conditions are:

\begin{iteMize}{$\bullet$}
\item It suffices to consider terms whose graphical representation
  does not contain any ``trivial cycles''. Trivial cycles are
  connected components of a diagram that do not contain any morphism
  variables. They can arise from the trace of an identity morphism. The
  restriction is without loss of generality because if $M$, $N$ have
  different numbers or types of trivial cycles, they can be easily
  separated in Hilbert spaces. Details are given in Lemma~\ref{lem-simple}
  below. We say that a diagram is {\em simple} if it contains no
  trivial cycles, and a term is simple if its associated diagram is
  simple.
\item We may assume that $M,N:I\ii I$, i.e., that both the domain and
  codomain of $M$ and $N$ are the tensor unit. Such terms are called
  {\em closed}. This simplification is justified in
  Lemma~\ref{lem-closed} below.
\item It suffices to consider terms $M,N$ in the language of {\em
    dagger traced monoidal} categories. Namely, by Joyal, Street, and
  Verity's ``Int''-construction {\cite{JSV96}}, every statement about
  dagger compact closed categories can be translated to an equivalent
  statement about dagger traced monoidal categories. Informally, this
  is done by eliminating occurrences of the $*$-operation: one
  replaces every morphism variable such as $f:A^*\otimes B\otimes
  C^*\ii D^*\otimes E$ by an equivalent new morphism variable such as
  $f':B\otimes D\ii A\otimes C\otimes E$ that does not use the
  $*$-operation. Details are given in Lemma~\ref{lem-tmc} below.
\end{iteMize}

\section{Informal outline of the result}\label{sec-informal}

The formal statement and proof of Theorem~\ref{thm-main} requires a
fair amount of notation. Nevertheless, the main idea is simple, and is
perhaps best illustrated in an example.  We thus start by giving an
informal explanation of the proof in this section.  The full technical
proof, including formal definitions of the graphical language and
isomorphism of diagrams, is given in Section~\ref{sec-technical}.

\subsection{Signatures, diagrams, and interpretations}

We assume given a set of {\em object variables}, denoted $A,B$ etc.,
and a set of {\em morphism variables}, denoted $f,g$ etc. A {\em sort}
$\vA$ is a finite sequence of object variables. We usually write
$A_1\otimes\ldots\otimes A_n$ for an $n$-element sequence, and $I$ for
the empty sequence. We assume that each morphism variable $f$ is
assigned two fixed sorts, called its {\em domain} $\vA$ and {\em
  codomain} $\vB$ respectively, and we write $f:\vA\ii\vB$. We further
require a fixed-point free involution $(-)\da$ on the set of morphism
variables, such that $f\da:\vB\ii\vA$ when $f:\vA\ii\vB$.

The collection of object variables and morphism variables, together
with the domain and codomain information and the dagger operation is
called a {\em signature} $\Sig$ of dagger monoidal categories.

Graphically, we represent a morphism variable
$f:A_1\otimes\ldots\otimes A_n\ii B_1\otimes\ldots\otimes B_m$ as a
box
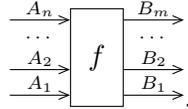
\[ \begin{tikzpicture}[scale=1,
    >=angle 45,
    box/.style={rectangle,draw,
    inner sep=1mm,minimum width=7mm,minimum height=12mm},
  ]
  \node[box,minimum height=13mm] (f) at (0,0) {$f$};
  \draw[<-] ([yshift=-.5cm]f.west) -- node[above,inner sep=1pt]{\tiny $A_1$} +(-.8cm,0);
  \draw[<-] ([yshift=-.16cm]f.west) -- node[above, inner sep=1pt]{\tiny $A_2$} +(-.8cm,0);
  \path ([yshift=.16cm]f.west) -- node[yshift=1mm]{\tiny $\cdots$} +(-.8cm,0);
  \draw[<-] ([yshift=.5cm]f.west) -- node[above,inner sep=1pt]{\tiny $A_n$} +(-.8cm,0);
  \draw[->] ([yshift=-.5cm]f.east) -- node[above,inner sep=1pt]{\tiny $B_1$} +(.8cm,0);
  \draw[->] ([yshift=-.16cm]f.east) -- node[above, inner sep=1pt]{\tiny $B_2$} +(.8cm,0);
  \path ([yshift=.16cm]f.east) -- node[yshift=1mm]{\tiny $\cdots$} +(.8cm,0);
  \draw[->] ([yshift=.5cm]f.east) -- node[above,inner sep=1pt]{\tiny $B_m$} +(.8cm,0);
\end{tikzpicture}.
\]
The wires on the left are called the {\em inputs} of $f$, and the
wires on the right are called its {\em outputs}. Note that each box is
labeled by a morphism variable, and each wire is labeled by an object
variable.

A (simple closed dagger traced symmetric monoidal) {\em diagram} over
a signature $\Sig$ consists of zero or more boxes of the above type,
all of whose wires have been connected in pairs, such that each
connection is between the output wire of some box and the input wire
of some (possibly the same, possibly another) box. Here is an example
of a diagram $N$ over the signature given by $f:B\ii A\otimes A$
and $g:A\otimes B\ii B\otimes A$.
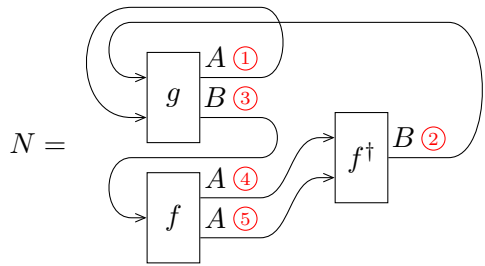
\[ \begin{tikzpicture}[scale=1,
  >=angle 45,
  box/.style={rectangle,draw,
    inner sep=1mm,minimum width=7mm,minimum height=12mm},
  input 1/.style={yshift=-.266cm},
  input 2/.style={yshift=.266cm},
  sp/.style={xshift=.8cm},
  ps/.style={xshift=-.2cm},
  obj/.style={above,pos=.25},
  number/.style={above,red,circle,draw,inner sep=1pt,outer sep=1mm,pos=.75}]
  \node[box] (f) at (0,0.2) {$f$};
  \node[box] (g) at (0,1.8) {$g$};
  \node[box] (fd) at (2.5,1) {$f^\dagger$};
  \node at (-1.8,1.2) {$N =$};
  \node[coordinate] (l1) at (0,2.8) {};
  \node[coordinate] (l2) at (0,3) {};
  \node[coordinate] (l3) at (0,1) {};
  \draw[->] ([input 1]f.east) -- node[obj]{$A$} node[number]{\tiny 5} ([input 1,sp]f.east) .. controls +(.4,0) and +(-.4,0).. ([input 1,ps]fd.west) -- ([input 1]fd.west);
  \draw[->] ([input 2]f.east) -- node[obj]{$A$} node[number]{\tiny 4} ([input 2,sp]f.east) .. controls +(.4,0) and +(-.4,0).. ([input 2,ps]fd.west) -- ([input 2]fd.west);
  \draw[->] (fd.east) -- node[obj]{$B$} node[number]{\tiny 2} ([sp]fd.east) .. controls +(.6,0) and +(.6,0).. ([sp]fd.east |- l1)
  -- ([ps]g.west |- l1) .. controls +(-.4,0) and +(-.4,0).. ([input 2,ps]g.west) -- ([input 2]g.west);
  \draw[->] ([input 2]g.east) -- node[obj]{$A$} node[number]{\tiny 1} ([input 2,sp]g.east) .. controls +(.4,0) and +(.4,0).. ([sp]g.east |- l2)
  -- ([ps]g.west |- l2) .. controls +(-.8,0) and +(-.8,0).. ([input 1,ps]g.west) -- ([input 1]g.west);
  \draw[->] ([input 1]g.east) -- node[obj]{$B$} node[number]{\tiny 3} ([input 1,sp]g.east) .. controls +(.3,0) and +(.3,0).. ([sp]g.east |- l3)
  -- ([ps]f.west |- l3) .. controls +(-.4,0) and +(-.4,0).. ([ps]f.west) -- (f.west);
\end{tikzpicture}
\]
In the illustration, we have numbered the wires $1$ to $5$ to aid the
exposition below; note that this numbering is not formally part of the
diagram.

An {\em interpretation} of a signature in finite-dimensional Hilbert
spaces consists of the following data: for each object variable $A$, a
chosen finite-dimensional Hilbert space $\semm{A}$, and for each
morphism variable $f:A_1\otimes\ldots\otimes A_n\ii
B_1\otimes\ldots\otimes B_m$, a chosen linear map
$\semm{f}:\semm{A_1}\otimes\ldots\otimes
\semm{A_n}\ii\semm{B_1}\otimes\ldots\otimes\semm{B_m}$, such that
$\semm{f\da}=\semm{f}\da$.

The {\em denotation} of a diagram $M$ under a given interpretation is
a scalar that is defined by the usual ``summation over internal
indices'' formula. For example, the denotation of the above diagram
$N$ is:
\begin{equation}\label{eqn-summation}
\semm{N} = \sum_{a_1,b_2,b_3,a_4,a_5}
\semm{g}\subsuper{a_1,b_2}{b_3,a_1}\cdot
\semm{f}\subsuper{b_3}{a_5,a_4}\cdot
\semm{f\da}\subsuper{a_5,a_4}{b_2}.
\end{equation}
Here $a_1, a_4, a_5$ range over some orthonormal basis of $\semm{A}$,
$b_2,b_3$ range over some orthonormal basis of $\semm{B}$, and
$\semm{f}\subsuper{b_3}{a_5,a_4}$ stands for the matrix entry
$\iprod{a_5\otimes a_4}{\semm{f}(b_3)}$. As is well-known, this
denotation is independent of the choice of orthonormal basis
{\cite{Geroch71,Penrose-Rindler87,Wald94}}.

\subsection{Proof sketch}\label{sec-proof-sketch}

By the reductions in Section~\ref{sec-reductions},
Theorem~\ref{thm-main} is a consequence of the following lemma:

\begin{lem}[Relative completeness]\label{lem-main}
  Let $M$ be a (simple closed dagger traced monoidal) diagram. Then
  there exists an interpretation $\semc{-}$ in finite dimensional
  Hilbert spaces, depending only on $M$, such that for all $N$,
  $\semc{N}=\semc{M}$ holds if and only if $N$ and $M$ are isomorphic
  diagrams.
\end{lem}

Clearly, the right-to-left implication is trivial, for if $N$ and $M$
are isomorphic diagrams, then $\semm{N}=\semm{M}$ holds under every
interpretation; their corresponding summation formulas differ at most
by a reordering of summands and factors. It is therefore the
left-to-right implication that must be proved.

The general proof of this lemma requires quite a bit of notation, as
well as more careful definitions than we have given above. A full
proof appears in Section~\ref{sec-technical} below. Here, we
illustrate the proof technique by means of an example.

Take the same signature as above, and suppose $M$ is the following diagram:
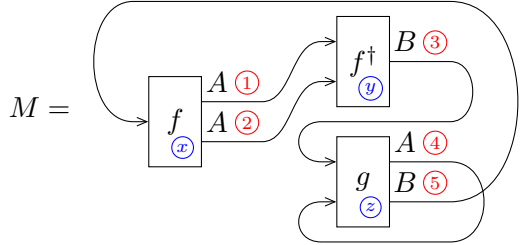
\[ \begin{tikzpicture}[scale=1,
  >=angle 45,
  box/.style={rectangle,draw,
    inner sep=1mm,minimum width=7mm,minimum height=12mm},
  input 1/.style={yshift=-.266cm},
  input 2/.style={yshift=.266cm},
  sp/.style={xshift=.8cm},
  ps/.style={xshift=-.2cm},
  obj/.style={above,pos=.25},
  number/.style={above,red,circle,draw,inner sep=1pt,outer sep=1mm,pos=.75},
  let/.style={below,blue,circle,draw,inner sep=1pt,outer sep=2mm,xshift=1mm}]
  \def\curve#1{.. controls +(#1,0) and +(#1,0)..}
  \node[box] (f) at (0,1) {$f$};
  \node[box] (g) at (2.5,.2) {$g$};
  \node[box] (fd) at (2.5,1.8) {$f^\dagger$};
  \node[let] at (f) {\tiny $x$};
  \node[let] at (fd) {\tiny $y$};
  \node[let] at (g) {\tiny $z$};
  \node at (-1.8,1.2) {$M=$};
  \node[coordinate] (l1) at (0,1) {};
  \node[coordinate] (l2) at (0,-.6) {};
  \node[coordinate] (l3) at (0,2.6) {};
  \draw[->] ([input 1]f.east) -- node[obj]{$A$} node[number]{\tiny 2} ([input 1,sp]f.east) .. controls +(.4,0) and +(-.4,0).. ([input 1,ps]fd.west) -- ([input 1]fd.west);
  \draw[->] ([input 2]f.east) -- node[obj]{$A$} node[number]{\tiny 1} ([input 2,sp]f.east) .. controls +(.4,0) and +(-.4,0).. ([input 2,ps]fd.west) -- ([input 2]fd.west);
  \draw[->] (fd.east) -- node[obj]{$B$} node[number]{\tiny 3} ([sp]fd.east) \curve{.4} ([sp]fd.east |- l1)
  -- ([ps]g.west |- l1) \curve{-.4} ([input 2,ps]g.west) -- ([input 2]g.west);
  \draw[->] ([input 2]g.east) -- node[obj]{$A$} node[number]{\tiny 4} ([input 2,sp]g.east) \curve{.6} ([sp]g.east |- l2)
  -- ([ps]g.west |- l2) \curve{-.4} ([input 1,ps]g.west) -- ([input 1]g.west);
  \draw[->] ([input 1]g.east) -- node[obj]{$B$} node[number]{\tiny 5} ([input 1,sp]g.east) \curve{1.2} ([sp]g.east |- l3)
  -- ([ps]f.west |- l3) \curve{-.7} ([ps]f.west) -- (f.west);
\end{tikzpicture}
\]
Again, we have numbered the wires from 1 to 5, and this time, we have
also numbered the boxes $x$, $y$, and $z$.

We must now construct the interpretation required by the Lemma. It is
given as follows. Define $\semc{A}$ to be a 3-dimensional Hilbert
space with orthonormal basis $\s{A_1,A_2,A_4}$. Define $\semc{B}$ to
be a 2-dimensional Hilbert space with orthonormal basis $\s{B_3,B_5}$.
Note that the names of the basis vectors have been chosen to suggest a
correspondence between basis vectors of $\semc{A}$ and wires labeled
$A$ in the diagram $M$, and similarly for $\semc{B}$.

Let $x$, $y$, and $z$ be three algebraically independent
transcendental complex numbers. This means that $x,y,z$ do not satisfy
any polynomial equation $p(x,y,z,\bar x, \bar y,\bar z)=0$ with
rational coefficients, unless $p\equiv 0$.

Define three linear maps $\semb{x}:\semc{B}\ii \semc{A}\otimes\semc{A}$,
$\semb{y}:\semc{A}\otimes\semc{A}\ii \semc{B}$,
and $\semb{z}:\semc{A}\otimes\semc{B}\ii \semc{B}\otimes\semc{A}$ as
follows. We give each map by its matrix representation in the chosen
basis.
\[\begin{array}{l}
  (\semb{x})\subsuper{i}{jk} = \left\{\begin{array}{l@{~~~}l}
      x&\mbox{if $i=B_5$, $j=A_2$, and $k=A_1$,}\\
      0&\mbox{otherwise,}
    \end{array}\right.\\
  (\semb{y})\subsuper{ij}{k} = \left\{\begin{array}{l@{~~~}l}
      y&\mbox{if $i=A_2$, $j=A_1$, and $k=B_3$,}\\
      0&\mbox{otherwise,}
    \end{array}\right.\\
  (\semb{z})\subsuper{ij}{kl} = \left\{\begin{array}{l@{~~~}l}
      z&\mbox{if $i=A_4$, $j=B_3$, $k=B_5$, and $l=A_4$,}\\
      0&\mbox{otherwise.}
    \end{array}\right.\\
\end{array}
\]
It is hopefully obvious how each of these linear functions is derived
from the diagram $M$: each matrix contains precisely one non-zero
entry, whose position is determined by the numbering of the input and
output wires of the corresponding box in $M$.

The interpretations of $f$ and $g$ are then defined as follows:
\[ \semc{f} = \semb{x} + \semb{y}\da,\quad\quad\semc{g} = \semb{z}.
\]
Note that we have taken the adjoint of the matrix $\semb{y}$, due to the
fact that the corresponding box was labeled $f\da$.
This finishes the definition of the interpretation $\semc{-}$.
It can be done analogously for any diagram $M$.

The reader may wonder why we didn't simply define $\semc{f} =
\semb{x}$ and $\semc{f\da} = \semb{y}$ independently of each
other. The reason is of course that an interpretation must satisfy
$\semc{f\da}=\semc{f}\da$.

To prove the condition of the Lemma, we first observe that the
interpretation $\semc{N}$ of any diagram $N$ is given by a summation
formula analogous to (\ref{eqn-summation}). Moreover, from the
definition of the interpretation $\semc{-}$, it immediately follows
that the scalar $\semc{N}$ can be (uniquely) expressed as a polynomial
$p(x,y,z,\bar x,\bar y,\bar z)$ with integer coefficients in the
variables $x,y,z$ and their complex conjugates. We note in passing
that this polynomial is homogeneous, and its degree is equal to the
number of boxes in $N$.

We claim that the coefficient of $p$ at $xyz$ is non-zero if and only
if $N$ is isomorphic to $M$. The proof is a direct calculation, using
(\ref{eqn-summation}) and the definition of $\semc{-}$.  Essentially,
any non-zero contribution to $xyz$ in the summation formula must come
from a choice of a basis vector $A_{\phi(w)}$ of $\semc{A}$ for each
wire $w$ labeled $A$ in $N$, and a choice of a basis vector
$B_{\phi(w)}$ of $\semc{B}$ for each wire $w$ labeled $B$ in $N$,
together with a bijection $\psi$ between the boxes of $N$ and the set
$\s{x,y,z}$; moreover, the contribution can only be non-zero if the
choice of basis vectors is ``compatible'' with the bijection $\psi$.
Compatibility amounts precisely to the requirement that the maps
$\phi$ and $\psi$ determine a graph isomorphism from $N$ to $M$. For
example, in the calculation of $\semc{N}$ according to equation
(\ref{eqn-summation}), the only non-zero contribution to $xyz$ in $p$
comes from the assignment $a_1\mapsto A_4$, $b_2\mapsto B_3$,
$b_3\mapsto B_5$, $a_4\mapsto A_1$, and $a_5\mapsto A_2$, which
corresponds exactly to the (in this case unique) isomorphism from $N$
to $M$.

In fact, we get a stronger result: the integer coefficient of $p$ at
$xyz$ is equal to the {\em number} of different isomorphisms between
$N$ and $M$ (usually 0 or 1, but it could be higher if $M$ has
non-trivial automorphisms).

\section{Technical development}\label{sec-technical}

\subsection{Signatures and diagrams}

\begin{defi}{\bf (Signature)}
  A {\em signature} of dagger monoidal categories is a quintuple
  $\Sig=\tuple{\Obj,\Mor,\dom,\cod,\dagger}$ consisting of:
  \begin{iteMize}{$\bullet$}
  \item a set $\Obj$ of {\em object variables}, denoted $A,B,C,\ldots$;
  \item a set $\Mor$ of {\em morphism variables}, denoted $f,g,h,\ldots$;
  \item functions $\dom,\cod:\Mor\ii\Obj^*$, called the {\em domain}
    and {\em codomain} functions, respectively, where $\Obj^*$ is the
    set of finite sequences of object variables;
  \item an operation $\dagger:\Mor\ii\Mor$, such that for all
    $f\in\Mor$, $f\dada=f$, $f\da\neq f$, $\dom f\da=\cod f$, and
    $\cod f\da=\dom f$.
  \end{iteMize}
\end{defi}

\noindent 
As before, we write a sequence of $n$ object variables as
$A_1\otimes\ldots\otimes A_n$, or as $\vA$, and we write $\abs{\vA}=n$
for the length of a sequence.  We write $f:\vA\ii \vB$ if $\dom f=\vA$
and $\cod f=\vB$.

\begin{defi}{\bf (Diagram)}
  A {\em (simple closed dagger traced symmetric monoidal) diagram}
  $M=\tuple{\WW{M},\BB{M},\labw{M},\labb{M},\thin{M},\thout{M}}$
  over a signature $\Sig$ consists of the following:
  \begin{iteMize}{$\bullet$}
  \item a finite set $\WW{M}$ of {\em wires};
  \item a finite set $\BB{M}$ of {\em boxes};
  \item a pair of labeling functions $\labw{M}:\WW{M}\ii\Obj$ and
    $\labb{M}:\BB{M}\ii\Mor$;
  \item a pair of bijections $\thin{M}:\Inputs{M}\ii\WW{M}$ and
    $\thout{M}:\Outputs{M}\ii\WW{M}$,
    where
    \[\begin{array}{lll}
      \Inputs{M} &=& \s{(i,b)\such b\in\BB{M}, n=\abs{\dom(\labb{M}(b))}, 1\leq i\leq n}, \\
      \Outputs{M} &=& \s{(b,j)\such b\in\BB{M}, m=\abs{\cod(\labb{M}(b))}, 1\leq j\leq m}. \\
    \end{array}
    \]
  \end{iteMize}
  Moreover, a diagram is required to satisfy the following typing conditions:
  \begin{iteMize}{$\bullet$}
  \item whenever $b\in\BB{M}$, $f=\labb{M}(b)$, $\vA=\dom f$, $(i,b)\in\Inputs{M}$, $\thin{M}(i,b)=w$, then $\labw{M}(w)=A_i$, and
  \item whenever $b\in\BB{M}$, $f=\labb{M}(b)$, $\vB=\cod f$, $(b,j)\in\Outputs{M}$, $\thout{M}(b,j)=w$, then $\labw{M}(w)=B_j$.
  \end{iteMize}
\end{defi}

\noindent 
Informally, $(i,b)$ represents the $i$th input of box $b$, $(b,j)$
represents the $j$th output of box $b$, and the bijections $\thin{M}$
and $\thout{M}$ determine which wires are attached to which inputs
and outputs, respectively. The labeling functions assign an object
variable to each wire and a morphism variable to each box, and the
typing conditions ensure that the sort of each wire matches the sort
of each box it is attached to.

\begin{defi}{\bf (Isomorphism)}
  An {\em isomorphism} of diagrams $N$, $M$ is given by a pair of
  bijections $\psiw:\WW{N}\ii\WW{M}$ and $\psib:\BB{N}\ii\BB{M}$,
  commuting with the labeling functions and with $\thin{}$ and
  $\thout{}$. Explicitly, this means that for all $w\in\WW{N}$
  and $b\in\BB{N}$, and all $i\leq \abs{\dom(\labb{N}(b))}$ and
  $j\leq\abs{\cod(\labb{N}(b))}$,
  \begin{eqnarray}
    \labw{M}(\psiw(w)) &=& \labw{N}(w), \label{eqn-iso1}\\
    \labb{M}(\psib(b)) &=& \labb{N}(b), \label{eqn-iso2}\\
    \thin{M}(i,\psib(b)) &=& \psiw(\thin{N}(i,b)), \label{eqn-iso3}\\
    \thout{M}(\psib(b),j) &=& \psiw(\thout{N}(b,j)). \label{eqn-iso4}
  \end{eqnarray}
\end{defi}

\begin{lem}\label{lem-redundant}
  In the definition of isomorphism, the condition that $\psiw$ is a
  bijection is redundant.
\end{lem}

\proof
  The bijection $\psib:\BB{N}\ii\BB{M}$ induces a bijection
  $\psibin:\Inputs{N}\ii\Inputs{M}$, defined by
  $\psibin(b,j)=(\psib(b),j)$. Equation (\ref{eqn-iso4}) is then
  equivalent to the commutativity of this diagram:
  \[ \xymatrix{
    \Inputs{N} \ar[r]^<>(.5){\thout{N}} \ar[d]_{\psibin} &
    \WW{N} \ar[d]^{\psiw} \\
    \Inputs{M} \ar[r]^<>(.5){\thout{M}} &
    \WW{M}.
    }
  \]
  Since the top, bottom, and left arrows are bijections, so is the
  right arrow.
\qed

\subsection{Interpretation in finite dimensional Hilbert spaces}

\begin{defi}{\bf (Interpretation)}
  Let $\Sig$ be a signature. An {\em interpretation} $\semm{-}$ of
  $\Sig$ in the category of finite dimensional Hilbert spaces assigns
  to each object variable $A\in\Obj$ a finite dimensional Hilbert
  space $\semm{A}$, and to each morphism variable $f:\vA\ii\vB$ a
  linear map $\semm{f}:\semm{A_1}\otimes \ldots\otimes \semm{A_n}\ii
  \semm{B_1}\otimes \ldots\otimes \semm{B_m}$, such that for all $f$,
  $\semm{f\da}=\semm{f}\da$.
\end{defi}

We sometimes write $\semm{\vA}$ for $\semm{A_1}\otimes \ldots\otimes
\semm{A_n}$.

\begin{defi}{\bf (Denotation)}
  Given a signature $\Sig$, a diagram $N$, and an interpretation
  $\semm{-}$. Fix an orthogonal basis $\Basis{A}$ for each space
  $\semm{A}$.  An {\em indexing} of $N$ is a function
  $\idx\in\prod_{w\in\WW{N}}\Basis{\labw{N}(w)}$, i.e., a choice of a
  basis element $\idx(w)\in\Basis{\labw{N}(w)}$ for every wire
  $w\in\WW{N}$. The set of indexings is written $\Idx{N}$. Then to
  each pair of an indexing $\idx$ and a box $b$, we assign a {\em
    matrix entry}
  \begin{equation}\label{eqn-bidx}
    b(\idx) = \semm{f}\subsuper
    {\idx(\thin{N}(1,b)),~\ldots,~\idx(\thin{N}(n,b)),}
    {\idx(\thout{N}(b,1)),~\ldots,~\idx(\thout{N}(b,m))}
  \end{equation}
  where $f=\labb{N}(b):A_1\otimes\ldots\otimes A_m\ii
  B_1\otimes\ldots\otimes B_m$.  As before, we have written
  $F\subsuper{y_1,\ldots,y_n}{x_1,\ldots,x_m}$ for the inner product
  $\iprod{x_1\otimes\ldots\otimes x_m}{F(y_1\otimes\ldots\otimes
    y_n)}$.  The {\em denotation} of $N$ is a scalar $\semm{N}$
  defined as follows:
  \begin{equation}\label{eqn-semmm}
    \semm{N}=\sum_{\idx\in\Idx{N}}\prod_{b\in\BB{N}} b(\idx).
  \end{equation}
\end{defi}

\begin{rem}
  The definition of $\semm{N}$ is independent of the chosen bases. In
  fact, the formula for $\semm{N}$ is just the usual formula for the
  summation over internal indices. It is the same as equation
  (\ref{eqn-summation}), expressed in the general context.
\end{rem}

\begin{rem}
  The graphical language can be interpreted in any dagger compact
  closed category {\cite{Sel05,Sel11}}. In the case of $\FinHilb$,
  the general interpretation coincides with the one given here.
\end{rem}

\subsection{The \texorpdfstring{$M$}{M}-interpretation}

We are now in a position to give a formal proof of relative
completeness (Lemma~\ref{lem-main}).  Fix a diagram $M$. We will
define a particular interpretation $\semc{-}$, called the {\em
  $M$-interpretation}, with the property that $\semc{N} = \semc{M}$ if
and only if $N$ and $M$ are isomorphic.

A family $\s{\xi_1,\ldots,\xi_k}$ of transcendental complex numbers is
{\em algebraically independent} if for every polynomial $p$ with
rational coefficients,
$p(\xi_1,\ldots,\xi_k,\bar\xi_1,\ldots,\bar\xi_k)=0$ implies $p\equiv
0$.  We choose a family of algebraically independent transcendental
numbers $\s{\xi_b\such b\in\BB{M}}$. This is possible because $\R$ is
a field extension of infinite transcendence degree over $\Q$
{\cite{Lang04}}, and because the complex numbers
$a_1+ib_1,\ldots,a_k+ib_k$ are algebraically independent in the above
sense if and only if $a_1,\ldots,a_k,b_1,\ldots,b_k$ are algebraically
independent real numbers.

For each object variable $A$, let $\WWW{M}{A}$ be the set of all wires
of $M$ that are labeled $A$, and for each morphism variable $f$, let
$\BBB{M}{f}$ be the set of boxes of $M$ that are labeled $f$. In
symbols,
\[\begin{array}{l}
  \WWW{M}{A}=\s{w\in \WW{M}\such \labw{M}(w)=A},\\
  \BBB{M}{f}=\s{b\in \BB{M}\such \labb{M}(b)=f}.\\
\end{array}
\]
Then the $M$-interpretation $\semc{-}$ is defined as follows. For each
$A$, let $\semc{A}$ be a Hilbert space with orthonormal basis
$\WWW{M}{A}$.  Suppose $f:\vA\ii\vB$ is a morphism variable, and
consider some $f$-labeled box $b\in\BBB{M}{f}$. We define a linear map
$\semb{b}:\semc{\vA}\ii\semc{\vB}$ by its matrix entries
\begin{equation}\label{eqn-sembb}
  (\semb{b})\subsuper
  {w_1,\ldots, w_n}
  {w'_1,\ldots, w'_m}
  = \left\{\begin{array}{l@{~~~}p{3in}}
      \xi_b & if $w_i=\thin{M}(i,b)$ and $w'_j=\thout{M}(b,j)$ for all $i,j$,\\
      0 & otherwise,
    \end{array}\right.
\end{equation}
where $w_i\in\WWW{M}{A_i}$ and $w'_j\in \WWW{M}{B_j}$ range over basis
vectors.  Finally, we define
\begin{equation}\label{eqn-semcf}
  \semc{f} = \sum_{b\in\BBB{M}{f}}\semb{b} +
  \sum_{b\in\BBB{M}{f\da}}\semb{b}\da.
\end{equation}

\subsection{Proof of relative completeness}

We must prove that the $M$-interpretation satisfies relative
completeness (Lemma~\ref{lem-main}).  First, we compute the
$M$-interpretation of any diagram $N$. By (\ref{eqn-semmm})
and (\ref{eqn-bidx}), we have
\[ \semc{N} = \sum_{\idx\in\Idx{N}}\prod_{b\in\BB{N}}
(\semc{\labb{N}(b)})\subsuper
{\idx(\thin{N}(1,b)),~\ldots,~\idx(\thin{N}(n,b)).}
{\idx(\thout{N}(b,1)),~\ldots,~\idx(\thout{N}(b,m))}
\]
Using (\ref{eqn-semcf}), it follows that
\[ \begin{array}{lll}
  \semc{N}
&=& \displaystyle\sum_{\idx\in\Idx{N}}\prod_{b\in\BB{N}}
\sum_{b'\in\BBB{M}{\labb{N}(b)}}
(\semb{b'})\subsuper
{\idx(\thin{N}(1,b)),~\ldots,~\idx(\thin{N}(n,b))}
{\idx(\thout{N}(b,1)),~\ldots,~\idx(\thout{N}(b,m))}
\\
&+& \displaystyle\sum_{\idx\in\Idx{N}}\prod_{b\in\BB{N}}
\sum_{b'\in\BBB{M}{\labb{N}(b)\da}}
(\semb{b'}\da)\subsuper
{\idx(\thin{N}(1,b)),~\ldots,~\idx(\thin{N}(n,b)).}
{\idx(\thout{N}(b,1)),~\ldots,~\idx(\thout{N}(b,m))}
\end{array}
\]
Now, using (\ref{eqn-sembb}) and the definition of $\dagger$, we
obtain the following explicit summation formula:
\[ \begin{array}{lll@{}l}
  \semc{N}
&=& \displaystyle\sum_{\idx\in\Idx{N}}\prod_{b\in\BB{N}}
\sum_{b'\in\BBB{M}{\labb{N}(b)}}
&\left\{
  \begin{array}{l@{~~~}l}
    \xi_{b'} &
    \parbox[t]{2.6in}{\raggedright
      if $\idx(\thin{N}(i,b))=\thin{M}(i,b')$ and
      $\idx(\thout{N}(b,j))=\thout{M}(b',j)$ for all $i,j$,}\\
    0 & \mbox{otherwise,}
  \end{array}
\right.\\
&+& \displaystyle\sum_{\idx\in\Idx{N}}\prod_{b\in\BB{N}}
\sum_{b'\in\BBB{M}{\labb{N}(b)\da}}
&\left\{
  \begin{array}{l@{~~~}l}
    \bar\xi_{b'} &
    \parbox[t]{2.6in}{\raggedright
      if $\idx(\thout{N}(b,i))=\thin{M}(i,b')$ and
      $\idx(\thin{N}(j,b))=\thout{M}(b',j)$ for all $i,j$,}\\
    0 & \mbox{otherwise.}
  \end{array}
\right.
\end{array}
\]
We note at this point that $\semc{N}$ can be (uniquely) written as a
polynomial with non-negative integer coefficients in the variables
$\s{\xi_b,\bar\xi_b\such b\in\BB{M}}$. Moreover, this polynomial is
homogeneous.

By definition, $b'\in\BBB{M}{\labb{N}(b)}$ if and only if
$\labb{M}(b')=\labb{N}(b)$, and $b'\in\BBB{M}{\labb{N}(b)\da}$ if and
only if $\labb{M}(b')=\labb{N}(b)\da$. The sets $\BBB{M}{\labb{N}(b)}$
and $\BBB{M}{\labb{N}(b)\da}$ are disjoint for each given $b$, since
$\dagger$ is fixed-point free. We can therefore rewrite the summation
as:
\[ \semc{N}
= \sum_{\idx\in\Idx{N}}\prod_{b\in\BB{N}}
\sum_{b'\in\BB{M}}
\left\{
  \begin{array}{l@{~~~}l}
    \xi_{b'} &
    \parbox[t]{2.6in}{\raggedright
      if $\labb{M}(b')=\labb{N}(b)$ and $\idx(\thin{N}(i,b))=\thin{M}(i,b')$ and
      $\idx(\thout{N}(b,j))=\thout{M}(b',j)$ for all $i,j$,}
    \\
    \bar\xi_{b'} &
    \parbox[t]{2.6in}{\raggedright
      if $\labb{M}(b')=\labb{N}(b)\da$ and $\idx(\thout{N}(b,i))=\thin{M}(i,b')$ and
      $\idx(\thin{N}(j,b))=\thout{M}(b',j)$ for all $i,j$,}
    \\
    0 & \mbox{otherwise.}
  \end{array}
\right.
\]
Finally, we use the distributive law to exchange the order of addition
and multiplication.
\begin{equation}\label{eqn-sumsumprod}
  \semc{N}
  = \sum_{\idx\in\Idx{N}}
  \sum_{\psii:\BB{N}\ii\BB{M}}
  \prod_{b\in\BB{N}}
  \left\{
    \begin{array}{l@{~~~}l}
      \xi_{\psii(b)} &
      \parbox[t]{2.6in}{\raggedright
        if $\labb{M}({\psii(b)})=\labb{N}(b)$ and $\idx(\thin{N}(i,b))=\thin{M}(i,{\psii(b)})$ and
        $\idx(\thout{N}(b,j))=\thout{M}({\psii(b)},j)$ for all $i,j$,}
      \\
      \bar\xi_{\psii(b)} &
      \parbox[t]{2.6in}{\raggedright
        if $\labb{M}({\psii(b)})=\labb{N}(b)\da$ and $\idx(\thout{N}(b,i))=\thin{M}(i,{\psii(b)})$ and
        $\idx(\thin{N}(j,b))=\thout{M}({\psii(b)},j)$ for all $i,j$,}
      \\
      0 & \mbox{otherwise.}
    \end{array}
  \right.
\end{equation}
Now consider a fixed $\idx\in\Idx{N}$ and fixed
$\psii:\BB{N}\ii\BB{M}$. We claim that the product
$\prod_{b\in\BB{N}}(\ldots)$ in (\ref{eqn-sumsumprod}) is equal to
$\prod_{b\in\BB{M}}\xi_b$ if and only if the pair of maps
$(\idx,\psii)$ forms an isomorphism of diagrams from $N$ to $M$.
Indeed, the product in question is equal to $\prod_{b\in\BB{M}}\xi_b$
if and only if $\psii$ is a bijection and the first side condition of
(\ref{eqn-sumsumprod}),
\begin{center}
  $\labb{M}({\psii(b)})=\labb{N}(b)$ and\\ $\idx(\thin{N}(i,b))=\thin{M}(i,{\psii(b)})$ and
  $\idx(\thout{N}(b,j))=\thout{M}({\psii(b)},j)$ for all $i,j$,
\end{center}
is satisfied for all $b\in\BB{N}$. This side condition amounts
precisely to conditions (\ref{eqn-iso2}), (\ref{eqn-iso3}), and
(\ref{eqn-iso4}) in the definition of isomorphism. Moreover, the
requirement that $\idx$ (viewed as a function from $\WW{N}$ to
$\WW{M}$) is an indexing is equivalent to condition (\ref{eqn-iso1}).
By Lemma~\ref{lem-redundant}, these conditions are necessary and
sufficient for the pair $(\idx,\psii)$ to be an isomorphism of
diagrams.

\proof[Proof of Lemma~\ref{lem-main}.]
  We have already noted that the right-to-left direction is trivial.
  For the left-to-right direction, note that we just showed that the
  polynomial $\semc{N}$ has a non-zero coefficient at
  $\prod_{b\in\BB{M}}\xi_b$ if and only if there exists an isomorphism
  between $N$ and $M$. In particular, $\semc{M}$ has a non-zero such
  coefficient. Therefore, if $\semc{N}=\semc{M}$, then $\semc{N}$ has
  a non-zero such coefficient, and it follows that $N$ and $M$ are
  isomorphic diagrams.
\qed

We note that the proof of Lemma~\ref{lem-main} yields a stronger
property:

\begin{cor}
  The coefficient of $\semc{N}$ at the monomial
  $\prod_{b\in\BB{M}}\xi_b$ is equal to the number of different
  isomorphisms from $N$ to $M$. \qed
\end{cor}

\subsection{Proof of the main result}

The following is an immediate consequence of Lemma~\ref{lem-main}.

\begin{prop}\label{prop-wlog}
  Let $M,N$ be two simple, closed terms in the language of dagger
  traced monoidal categories. Suppose that $\semm{M}=\semm{N}$ for
  every possible interpretation in finite dimensional Hilbert
  spaces. Then $M=N$ holds in the graphical language.
\end{prop}

The proof of Theorem~\ref{thm-main} then consists of removing the
conditions ``simple'', ``closed'', and ``in the language of dagger
traced monoidal categories'', using the reductions outlined in
Section~\ref{sec-reductions}.  We give the details here.

\begin{lem}\label{lem-simple}
  The condition ``simple'' can be removed from
  Proposition~\ref{prop-wlog}.
\end{lem}

\proof
  Let $M,N:A\ii B$ be two terms, not necessarily simple. It suffices
  to show that if there is some object variable $A$ such that $M$ and
  $N$ have a different number of trivial cycles of type $A$, then
  there exists an interpretation such that $\semm{M}\neq\semm{N}$. Let
  $\semm{A}$ be a 2-dimensional Hilbert space with basis
  $\s{a,b}$. For any $B\neq A$, let $\semm{B}$ be a 1-dimensional
  Hilbert space with basis $\s{c}$. We call the basis vectors $a$ and
  $c$ {\em distinguished}. For a morphism variable $f:\vA\ii\vB$,
  define
  \[ \semm{f}\subsuper{w_1,\ldots,w_n}{w'_1,\ldots,w'_m} =
  \left\{\begin{array}{l@{~~~}p{3.4in}}
      1 & if $w_1,\ldots,w_n,w'_1,\ldots,w'_m$ are all distinguished,\\
      0 & otherwise,
    \end{array}\right.
  \]
  It is then easy to see that any {\em simple} closed term $P$
  satisfies $\semm{P}=1$. Namely, in (\ref{eqn-semmm}),
  $\prod_{b\in\BB{N}} b(\idx)$ is equal to $1$ precisely for the
  unique indexing $\idx\in\Idx{N}$ that picks only distinguished basis
  vectors, and $0$ otherwise. Also, the interpretation of a trivial
  cycle of type $A$ is $\Tr(\id_A) = \dim\semm{A} = 2$, whereas the
  interpretation of any other trivial cycle is
  $\dim\semm{B}=1$. Therefore, if $M$ and $N$ have $k$ and $l$ trivial
  cycles of type $A$, respectively, then $\semm{M}=2^k$ and
  $\semm{N}=2^l$, which implies $\semm{M}\neq\semm{N}$ if $k\neq l$.
  \qed

\begin{lem}\label{lem-closed}
  The condition ``closed'' can be removed from
  Proposition~\ref{prop-wlog}.
\end{lem}

\proof
  Suppose Proposition~\ref{prop-wlog} holds for all closed terms.  Let
  $M,N:A\ii B$ be two terms that are not necessarily closed, such that
  every possible interpretation satisfies $\semm{M}=\semm{N}$. We can extend
  the language with two new morphism variables $f:I\ii A$ and $g:B\ii
  I$, and apply the proposition to the terms $M'=g\cp M\cp f$ and
  $N'=g\cp N\cp f$. It follows that the diagrams for $g\cp M\cp f$ and
  $g\cp N\cp f$ are isomorphic in the graphical language. But since
  $f,g$ are new symbols that do not occur in $M$ and $N$, this implies
  that $M$ and $N$ are isomorphic as well. 
\qed

\begin{lem}\label{lem-tmc}
  Proposition~\ref{prop-wlog} holds for terms in the language of
  dagger compact closed categories, over any signature of dagger
  compact closed categories.
\end{lem}

\proof Let $\Sig$ be a signature of dagger compact closed
categories. This means that the domain and codomain of each morphism
variable is a sequence of object variables and their duals, for
example $f:A^*\otimes B\otimes C^*\ii D^*\otimes E$. Let $\CCSig$ be
the free dagger compact closed category over $\Sig$.  The first
observation is that there exists a signature $\Sig'$ of dagger {\em
  monoidal} categories such that $\CCSigp\iso\CCSig$ are isomorphic
categories. Indeed, $\Sig'$ is obtained from $\Sig$ by replacing each
$A^*$ in the domain of a morphism variable by $A$ in its codomain, and
vice versa; for example, the above $f$ in $\Sig$ will be replaced by
$f':B\otimes D\ii A\otimes C\otimes E$ in $\Sig'$. Then
$\CCSigp\iso\CCSig$ because $f$ and $f'$ are interdefinable in any
compact closed category.

  What we must show is that for all $f\neq g:A\ii B$ in $\CCSig$,
  there exists some dagger compact closed functor
  $F:\CCSig\ii\FinHilb$ such that $F(f)\neq F(g)$. Equivalently, we
  have to show that there exists a faithful dagger compact closed
  functor $\bar F:\CCSig\ii\FinHilb^X$ into some discrete power
  $\FinHilb^X = \prod_{i\in X}\FinHilb$ of the category of finite
  dimensional Hilbert spaces.

  Let $\TMCSigp$ be the free dagger traced symmetric monoidal category
  over $\Sig'$. It is an easy exercise to prove that Joyal, Street,
  and Verity's ``Int''-construction {\cite{JSV96}} freely embeds any
  dagger traced symmetric monoidal category $\Dd$ in a dagger compact
  closed category $\Int(\Dd)$, and moreover, that any faithful dagger
  traced monoidal functor $G:\Dd\ii\Cc$ into a dagger compact closed
  category extends to a faithful dagger compact closed functor $\hat
  G:\Int(\Dd)\ii\Cc$. Applying this to the situation where
  $\Dd=\TMCSig$ and $\Cc=\CCSigp$, and using the respective universal
  properties of $\CCSigp$ and of the Int-construction, we obtain an
  equivalence of dagger compact closed categories
  $\CCSigp\simeq\Int(\TMCSig)$.

  By hypothesis (i.e., Proposition~\ref{prop-wlog}, with the
  conditions ``simple'' and ``closed'' already removed), there exists
  a faithful traced monoidal functor $H:\TMCSigp\ii\FinHilb^X$ for
  some $X$. Then the composition 
  \[ \bar F :=
  \CCSig\catarrow{\iso}\CCSigp\catarrow{\simeq}\Int(\TMCSig)\catarrow{\hat
    H}\FinHilb^X
  \] is the desired faithful dagger compact closed
  functor. 
\qed

\section{Generalizations}\label{sec-generalizations}

\paragraph{Other rings and fields}

The result of this paper (Theorem~\ref{thm-main}) can be adapted to
other fields besides the complex numbers.  It is true for any field
$k$ of characteristic $0$ with a non-trivial involutive automorphism
$x\mapsto \bar x$. (Non-trivial means that for some $x$, $\bar x\neq
x$).

The only special property of $\C$ that was used in the proof, and
which may not hold in a general field $k$, was the existence of
transcendentals.  This problem is easily solved by first considering
the field of fractions $k(x_1,\bar x_1\ldots,x_n,\bar x_n)$, where the
required transcendentals have been added freely. The proof of
Lemma~\ref{lem-main} then proceeds without change. Finally, once an
interpretation over $k(x_1,\bar x_1\ldots,x_n,\bar x_n)$ has been
found such that $\semm{M}\neq\semm{N}$, we use the fact that in a
field of characteristic 0, any non-zero polynomial has a
non-root. Thus we can instantiate $x_1,\ldots,x_n$ to specific
elements of $k$ while preserving the inequality
$\semm{M}\neq\semm{N}$. Note that therefore, Theorem~\ref{thm-main}
holds for the given field $k$; however, Lemma~\ref{lem-main} only holds for
$k(x_1,\bar x_1\ldots,x_n,\bar x_n)$.

Moreover, the elements $x_1,\ldots,x_n$ of the preceding paragraph can
always be instantiated to integers; therefore, the results also hold
if one replaces $k$ by the ring of Gaussian integers $k=\Z[i]$. In
elementary terms: if a certain equation fails to hold in dagger
compact closed categories, then one can always find a counterexample
in matrices with entries of the form $a+bi$, where
$a,b\in\Z$. Moreover, $i=\sqrt{-1}$ can be replaced by $\sqrt{d}$ for
any non-square integer $d$.

\paragraph{Non-equational properties}

The completeness result of this paper applies to properties
expressible as equations. One may ask whether it can be generalized to
other classes of properties, such as implications, Horn clauses, or
more general logical formulas. Unfortunately, this is not the case;
there exist (universally quantified) implications that hold in
$\FinHilb$, but fail to hold in arbitrary dagger compact closed
categories. One example of such an implication is
\[ \forall f,g:A\to B\ssep (f\da ff\da f = g\da gg\da g \ssep\imp\ssep f\da f = g\da g).
\]
This is true in $\FinHilb$, because each hermitian positive operator
has a unique hermitian positive square root. But it fails in general
dagger compact closed categories. Perhaps the simplest counterexample
is the ring $\Z_5$, regarded as a one-object dagger compact closed
category with composition and tensor of morphisms given by
ring multiplication, and dagger given by the identity operation. Taking
$f=1$ and $g=2$, the premise is satisfied and the conclusion is not.

\paragraph{Bounded dimension}

The interpretation $\semc{-}$ from the proof of Lemma~\ref{lem-main}
uses Hilbert spaces of unbounded finite dimension. One may ask whether
Theorem~\ref{thm-main} remains true if the interpretation of object
variables is restricted to Hilbert spaces of some fixed dimension
$n$. This is known to be false when $n=2$. Here is a counterexample
due to Bob Par\'e: the equation $\tr(AABBAB)=\tr(AABABB)$ holds for
all $2\times 2$-matrices, but does not hold in the graphical
language. Indeed, by the Cayley-Hamilton theorem, $A^2=\mu A+\nu I$
for some scalars $\mu,\nu$. Therefore
\[\begin{array}{l}
  \tr(AABBAB) = \mu\tr(ABBAB) + \nu\tr(BBAB), \\
  \tr(AABABB) = \mu\tr(ABABB) + \nu\tr(BABB),
\end{array}
\]
and the right-hand-sides are equal by cyclicity of trace. It is not
currently known to the author whether Theorem~\ref{thm-main} is true
when restricted to spaces of dimension~3.

\section*{Acknowledgments}

Thanks to Gordon Plotkin for telling me about this problem and for
discussing its solution. Thanks to Bob Par\'e for the counterexample in
dimension 2.

\bibliographystyle{abbrv}
\bibliography{finhilb}

\end{document}